\documentclass[twocolumn,preprintnumbers]{revtex4}
\usepackage{graphicx,color}
\usepackage{dcolumn}
\usepackage{bm}

\newcommand{\be}{\begin{equation}}
\newcommand{\ee}{\end{equation}}
\newcommand{\bea}{\begin{eqnarray}}
\newcommand{\bc}{\begin{center}}
\newcommand{\ec}{\end{center}}
\newcommand{\eea}{\end{eqnarray}}

\draft
\input{tcilatex}
\begin{document}

\title{Genaral Pi Function and the Pi number }
\author{Hakan Ciftci \\
Gazi \"{U}niversitesi, Fen-Edebiyat Fak\"{u}ltesi, \\
Fizik B\"{o}l\"{u}m\"{u}, 06500, Teknikokullar, Ankara/TURKEY\\
\ }

\begin{abstract}
In this short note, we have defined a new "nested square root"
function which generates usual Pi number for $x=2$. We have given
some useful identities and asymptotic formulas of the Pi-function.
\end{abstract}

\maketitle

\section{Introduction}

Studies on new methods calculating the usual \textbf{Pi-number} have
been done so far, It is possible to find many works on this issue
[1-10]. In this work our main interest is not getting the Pi-number,
but obtaining more general Pi-function which gives us usual
Pi-number for special case of the general function. This kind of
definition will be helpful to get different representation of the
number Pi. To be able to get this function, we will follow
Vi\`{e}te while we are obtaining the general $\Pi (x)$ function. In 1543, Vi%
\`{e}te showed that usual Pi-number could be obtained as the following
nested square root%
\begin{equation}
{\pi =\lim_{i\rightarrow \infty }(4)}^{\frac{i+1}{2}}\sqrt{2-h_{i}}
\end{equation}%
where
\begin{equation}
{h}_{i}=\sqrt{2+h_{i-1}}\text{, with }h_{0}=0
\end{equation}

One can obtain this formula using the following iterative procedure

\begin{enumerate}
\item First, suppose that we have a circle and its radius is $\sqrt{2}$. We
can draw a square inside this circle and let the corners of that square
touch the circle. Thus, the length of the one side of the square is $2$. Now
if we calculate the ratio of the circumference of the square to the diameter
of the circle, we get%
\begin{equation}
{\pi \approx }\frac{8}{2\sqrt{2}}=2\sqrt{2}
\end{equation}

\item Second, we can draw octagon inside the circle and let the corners of
the octagon touch the circle. The radius of the circle is still $\sqrt{2}$.
Thus we can easily calculate the length of the one side of the octagon as
below
\begin{equation}
{x=}\sqrt{2}\sqrt{2-\sqrt{2}}
\end{equation}%
an later we can again calculate the ratio of the circumference of the
octagon to the diameter of the circle as below%
\begin{equation}
{\pi \approx }4\sqrt{2-\sqrt{2}}
\end{equation}
\end{enumerate}

If these calculations are performed, one can easily find Eq.(1)

\section{Generalization}

In this section, we would like to write a general equation which covers
Eq.(1). For this purpose, we will make the following ansatz,
\begin{equation}
\Pi _{i}(x){=(2x)}^{\frac{i+1}{2}}\sqrt{x-h_{i}}\text{, }x\in
(-\infty,\infty )
\end{equation}%
where
\begin{equation}
{h}_{i}=\sqrt{x(x-1)+h_{i-1}}\text{, with }h_{0}=0
\end{equation}%
This $\Pi _{i}(x)$ function approaches to a certain value as $i$ goes to
infinity if $x\in \lbrack 1^{+},\infty )$. That is to say:%
\begin{equation}
\lim_{i\rightarrow \infty }\Pi _{i}(x){=\pi }_{x}
\end{equation}%
It is easy to see that when we take $x=2$ in Eq(6) and (7), we get Eq.(1).
Eq.(6) generates interesting numbers for us. We should point out that if $%
z=x+iy$ $(x\neq 0)$, then we can get that $\Pi _{i}(z)$ approaches
to same positive and negative constant complex numbers. In Table-1
we give some numerical example for $\Pi _{i}(x)$ and in the Figure
the behavior of the \textit{Pi-function} has been shown.
\begin{figure}[tbp]
\centering
\includegraphics[width=4.5cm,clip=true]{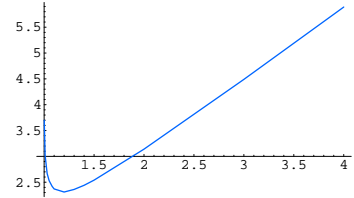}
\caption{Behaviour of the $\Pi (x)$}
\label{fig:fig1}
\end{figure}
To use this $\Pi -$function may be useful to produce different
calculation ways for usual Pi-number. We will give some interesting
properties of the \textit{Pi-function} and later using these
identities, one can to obtain different representation of the
pi-number.

\subsection{Some properties of $\Pi _{i}(x)$}

It is easy to show that
\begin{equation}
\lim_{x\rightarrow \infty }\Pi _{i}(x){=}\sqrt{2}x\left( 1+\frac{1}{8x}%
\right)
\end{equation}%
and we have also determined that $\Pi _{i}(x)$ function has a minimum value
at $x\approx 1.19005.$ At this value,
\begin{equation}
\Pi _{i}(1.19005){\approx 2.31383}
\end{equation}%
Another interesting properties can be written as below%
\begin{equation}
\lim_{x\rightarrow \infty }\sqrt{\frac{x-h_{i}}{x-h_{i+1}}}{=}\sqrt{2x}
\end{equation}%
where $h_{j}$ is defined as Eq.(2). In addition to above identities, we can
obtain some new formulas to find number-Pi, for this purpose, one can
calculate $\Pi _{i+1}(x)$ as below%
\begin{equation}
\Pi _{i+1}^{2}(x){=2x}\Pi _{i}^{2}(x)+(2x)^{i+2}(h_{i}-h_{i+1})
\end{equation}%
From this formula, one can easily get that when $i$ goes to infinity, $\Pi
_{i+1}$ and $\Pi _{i}$ approaches to the same $\pi _{x}$ values given in
Eq.(8). Thus, one can write the following formula while $i$ goes to infinity%
\begin{equation}
\pi _{x}=\left( 2x\right) ^{\frac{i+2}{2}}\sqrt{\frac{h_{i+1}-h_{i}}{2x-1}}
\end{equation}%
If we calculate asymptotic behavior of Eq.(13) at large $x$ values one can
obtain that
\begin{equation}
\pi _{x}=x\sqrt{\frac{4x-1}{2x-1}}
\end{equation}%
Additionally, when we look at the $h_{i}(x)$ function, we see that
\begin{equation}
h_{i}(x)=h_{i}(1-x)
\end{equation}%
from this equation, one can calculate that
\begin{equation}
\Pi _{i}^{2}(1-x){=(}\frac{1}{x}-1{)}^{i+1}\Pi
_{i}^{2}(x)-2^{i+1}(1-x)^{i+1}(2x-1)
\end{equation}

Table.1 Calculation of Pi function for different values of $x$. Iteration
number ($i$) is taken as 50

\begin{tabular}{|l|l|}
\hline
$x$ & $\Pi _{i}(x)$ \\ \hline
$1.001$ & $3.7033451$ \\ \hline
$1.5$ & $2.5351046$ \\ \hline
$\mathbf{2}$ & $\mathbf{3.1415927}$ \\ \hline
$2.5$ & $3.8084662$ \\ \hline
$3$ & $4.4937674$ \\ \hline
$4$ & $5.8848462$ \\ \hline
$5$ & $7.2869301$ \\ \hline
$7$ & $10.102809$ \\ \hline
$8$ & $11.513355$ \\ \hline
$20$ & $28.469656$ \\ \hline
\end{tabular}

\bigskip

\section{References}
\begin{enumerate}

\item T. J. Osler , FIBONACCI QUARTERLY 45, (2007) 202.

\item P. J. Humphries, BULLETIN OF THE KOREAN MATHEMATICAL SOCIETY 44,
(2007) 331

\item J. Johnson, T. Richmond, RAMANUJAN JOURNAL 15, (2008) 259

\item A. Gee, M. Honsbeek, RAMANUJAN JOURNAL 11 (2006) 267

\item A. Levin, RAMANUJAN JOURNAL 10, (2005) 305

\item K. S. Rao, G. V. Berghe, JOURNAL OF COMPUTATIONAL AND APPLIED
MATHEMATICS 173 (2005) 371

\item J. Blomer, ALGORITHMICA 28, (2000) 2

\item N. N. Osipov, PROGRAMMING AND COMPUTER SOFTWARE 23, (1997) 142

\item S.Landau, SIAM JOURNAL ON COMPUTING 21, (1992) 85

\item J. M. Borwein, G. Debarra, AMERICAN MATHEMATICAL MONTHLY 98, (1991) 735
\end{enumerate}

\end{document}